\documentclass[12pt]{article}
\usepackage{amssymb, amsfonts, amsmath}
\usepackage{hyperref}
\usepackage{xcolor}
\hypersetup{
	colorlinks,
	linkcolor={red!60!black},
	citecolor={green!60!black}
}
\newcommand{\qed}{\hskip 5mm \rule{2.5mm}{2.5mm}\vskip 10pt}
\newcommand{\R}{{\mathbb R}}

\newcommand{\E}{{\mathbb E}}

\newcommand{\proof}[1]{\noindent{\em Proof:\ } #1 \qed}

\newcommand{\abs}[1]{ \left| #1 \right|}

\newcommand{\func}[2]{#1 \left( #2 \right)}

\newcommand{\parent}[1]{ \left( #1 \right)}

\newcommand{\braces}[1]{ \left\{ #1 \right\}}
\newcommand{\setbuilder}[2]{ \left\{ #1 \mid #2 \right\}}

\newcommand{\3}[1]{\mathcal{#1}}

\numberwithin{equation}{section}

\begin{document}
\newtheorem{theorem}{Theorem}[section]
\newtheorem{definition}[theorem]{Definition}
\newtheorem{lemma}[theorem]{Lemma}
\newtheorem{note}[theorem]{Note}
\newtheorem{corollary}[theorem]{Corollary}
\newtheorem{proposition}[theorem]{Proposition}
\renewcommand{\theequation}{\arabic{section}.\arabic{equation}}
\newcommand{\newsection}[1]{\setcounter{equation}{0} \section{#1}}
%%%%%%%%%%%% title %%%%%%%%%%%%%%%%%%%%%%%%%%%%%%%%
\title{Ergodicity in Riesz spaces
}
%%%%%%%%%%%%%%%%%%%%%%%%%%%%%%%%%%%%%%%%%%%%%%%%
\author{
 Jonathan Homann\\
 Wen-Chi Kuo\\
 Bruce A. Watson \\ \\
 ${}^\sharp$ School of Mathematics\\
 University of the Witwatersrand\\
 Private Bag 3, P O WITS 2050, South Africa }

\maketitle
\abstract{
The ergodic theorems of Hopf, Wiener and Birkhoff were extended to the context of Riesz spaces with a weak order unit and conditional expectation operator by Kuo, Labuschagne and Watson in  
[\textit{Ergodic Theory and the Strong Law of Large Numbers on Riesz Spaces.} 
{{Journal of Mathematical Analysis and Applications}}, 
{\textbf{325}}, 
{(2007), 422--437}.]. However, the precise concept of what constitutes ergodicity in Riesz spaces was not considered. In this short paper we fill in this omission and give some explanations of the choices made. In addition, we consider the interplay between mixing and ergodicity in the Riesz space setting.
}
\parindent=0cm
\parskip=0.5cm
%%%%%%%%%%%%%%%%%%%%%%%%%%%%%%%%%%%%%%%%%%%%%%
\section{Introduction} \label{s: introduction}

Ergodic theory seeks to study the long-term behaviour of a dynamical system. One typically works in a probability space, $\parent{\Omega,\3A,\mu}$, with a transformation $\tau \colon \Omega \rightarrow \Omega$. Here the quadruple $\parent{\Omega,\3A,\mu,\tau}$ is called a measure preserving system if $\func{\mu}{{\tau}^{-1}(A)}=\func{\mu}{A}$ for each $A \in \3A$.  A set $A \in \3A$ is called $\tau$-invariant if $\func{{\tau}^{-1}}{A}\subset A$ and the
measure preserving system $\parent{\Omega,\3A,\mu,\tau}$ is said to be ergodic if $\func{\mu}{A} \in \braces{0,1}$ for each $\tau$-invariant $A$, see Petersen \cite{ergodic book}. An equivalent definition of ergodicity is that $f\circ \tau = f$ for measurable $f$ only if $f$ is a.e. constant, see \cite[page 15]{shields}. It is this definition that leads naturally to a conditional version and, from there, to a Riesz space version.

In \cite{ergodic paper}, Kuo, Labuschagne and Watson generalised the ergodic theorems of Birkhoff, Wiener, Hopf and Garsia to the setting of Riesz spaces. They also generalised the Kolmogorov 0-1 law to Riesz spaces and were thus able to give a strong law of large numbers for conditionally independent sequences in Riesz spaces. However, they omitted to consider what constituted ergodicity in a Riesz space. In this paper we fill this gap and, using the works of Kuo, Labuschagne and Watson in \cite{ergodic paper, khw-1-A},  are able to relate ergodicity to weak mixing in the Riesz space setting. Further to this, we apply these concepts to sequences with conditionally independent shifts and show the ergodicity of such processes.

This work continues the studies of mixing processes, \cite{khw-1-A, krw}, and of stochastic processes, \cite{azouzi, jensen paper, expectation paper, ergodic paper, KRodW2, stoica}, in the Riesz space setting.

This work is dedicated to the memory of our friend, collaborator and colleague  Professor Coenraad Labuschagne, with whom our work on stochastic processes in Riesz spaces began.
%%%%%%%%%%%%%%%%%%%%%%%%%%%%%%%%%%%%%%%%%%%%%%%%%%%%%%%
\section{Preliminaries} \label{preliminary section}
The reader is referred to Aliprantis and Border \cite[pages 263-300]{riesz book}, Fremlin \cite[Chapter 35, pages 219-274]{fremlin}, Meyer-Nieberg \cite{MN-BL}, and Zaanen \cite{zaanen} for background in Riesz spaces.

We recall from \cite{expectation paper} the definition of a conditional expectation operator on a Riesz space.

\begin{definition}
Let $E$ be a Riesz space with weak order unit. A positive order continuous projection $T \colon E \rightarrow E$, with range, $\func{R}{T}$, a Dedekind complete Riesz subspace of $E$, is called a conditional expectation operator on $E$ if $Te$ is a weak order unit of $E$ for each weak order unit $e$ of $E$.
\end{definition}

The Riesz space analogue of a measure preserving system is introduced in the following definition.

\begin{definition} \label{system definition}
If $E$ is a Dedekind complete Riesz space with weak order unit, say, $e$, $T$ is a conditional expectation operator on $E$ with $Te=e$ and $S$ is an order continuous Riesz homomorphism on $E$ with $Se=e$ and, further, $TSPe=TPe$, for all band projections $P$ on $E$, then $\parent{E,T,S,e}$, is called a conditional expectation preserving system.
\end{definition}

Due to Freudenthal's Spectral Theorem (\cite[Theorem 33.2]{zaanen}), the condition $TSPe=TPe$ for each band projection $P$ on $E$ in the above definition is equivalent to $TSf=Tf$ for all $f \in E$.

Let  $(E,T,S,e)$ be a conditional expectation preserving system.
For $f \in E$ and $n \in \mathbb{N}$ we denote
	\begin{subequations} \label{eqn def}
	\begin{align}
	S_nf &:= \frac{1}{n}\sum_{k=0}^{n-1} S^kf,
	\label{s_n eqn def}\\
	L_Sf &:= \lim_{n \rightarrow \infty} S_n f,	
	\label{fhat eqn def}
	\end{align}
	\end{subequations}
where the above limit is the order limit, if it exits.
We say that  $f \in E$ is $S$-invariant if $Sf=f$. 
The set of all $S$-invariant $f \in E$ will be denoted 
$\3I_S := \setbuilder{f \in E}{Sf=f}$.
The set of $f \in E$ for which $L_S f$ exists will be denoted $\3E_S$ and thus $L_S:\3E_S\to E$.

\begin{lemma}
\label{independent exist lemma}
Let $\parent{E,T,S,e}$ be a conditional expectation preserving system and define $\3I_S$ and $\3E_S$ as before, then $\displaystyle \3I_S \subset \3E_S$ and $L_Sf=f$, for all $f\in \3I_S$, so $\3I_S\subset R(L_S)$.
\end{lemma}

\proof{
If $f \in \3I_S$, then $S^kf=f$, for all $k \in \mathbb{N}_0$. Thus  $\displaystyle S_nf=\frac{1}{n}\sum_{k=0}^{n-1}f=f$, for all $n \in \mathbb{N}$, giving $S_n f\to f$, in order, as $n\to\infty$. Hence $f\in \3E_S$ and
$L_Sf=f$, thus $\3I_S\subset \3E_S$.
}

We recall Birkhoff's bounded ergodic theorem from \cite[Theorem 3.7]{ergodic paper}.

\begin{theorem}[Birkhoff's (Bounded) Ergodic Theorem] \label{bounded riesz birkhoff theorem}

Let  $(E,T,S,e)$ be a conditional expectation preserving system.
For $f\in E$, the sequence $\displaystyle {\parent{S_n f}}_{n \in \2N}$ is order bounded in $E$ if and only if $f\in \3E_S$.
For each $f\in\3E_S$ we have $L_Sf=SL_Sf$ and $TL_Sf=Tf$. If $E=\3E_S$ then $L_S$ is a conditional expectation operator on $E$ with $L_Se=e$.
\end{theorem}

\begin{note}
If we restict our attention to $E_e$, then $\displaystyle {\parent{S_n f}}_{n \in \2N}$ is order bounded in $E_e$ for each $f\in E_e$, so $E_e \subset \3E_S$. Furthermore, $L_Sf \in E_e$, giving that $L_S|_{E_e}$ is a conditional expectation on $E_e$.
\end{note}

$E$ is said to be universally complete with respect to $T$ ($T$-universally complete), if, for each increasing net ${\parent{f_{\alpha}}}$ in $E_+$ with $\parent{Tf_{\alpha}}$ order bounded, we have that ${\parent{f_{\alpha}}}$ is order convergent in $E$.

We recall Birkhoff's ergodic theorem for a $T$-universally complete Riesz space from \cite[Theorem 3.9]{ergodic paper}.

\begin{theorem}[Birkhoff's (Complete) Ergodic Theorem]
\label{complete riesz birkhoff theorem}
Let $(E,T,S,e)$ be a conditional expectation preserving system with $E$ $T$-universally complete then $E=\3E_S$ and hence $L_S=SL_S$. In addition, $TL_S=T$ and $L_S$ is a conditional expectation operator on $E$.
\end{theorem}

From \cite[Lemma 2.1]{khw-1-A} we have.

\begin{lemma} \label{absolute-cesaro}
Let $\parent{f_n}_{n \in \2N_0}$ be a sequence in $E$, a Dedekind complete Riesz space. If  
$\displaystyle{\sum_{k=0}^{\infty}|f_k|}$ is order convergent in $E$ then  
$\displaystyle{\sum_{k=0}^{\infty}f_k}$ is order convergent in $E$.
\end{lemma}

From the above lemma and \cite[Lemma 2.1]{KRodW2} we have the following.

\begin{theorem} \label{convergent sequence-absolute cesaro sum}
Let $E$ be a Dedekind complete Riesz space and $\parent{f_n}_{n \in \2N_0}$ a sequence in $E$ with $f_n \to 0$, in order, as $n \to \infty$, then
$\displaystyle \frac{1}{n} \sum_{k=0}^{n-1} \abs{f_k}\to 0,$
in order, as $n \to \infty$.
\end{theorem}

\begin{corollary} \label{convergence implies cesaro convergence}
Let $E$ be a Dedekind complete Riesz space and $\parent{f_n}_{n \in \2N_0}$ a sequence in $E$ with order limit $f$, then
$\displaystyle
\frac{1}{n} \sum_{k=0}^{n-1}f_k
\to f,$ in order, as $n \to \infty$.
\end{corollary}

%%%%%%%%%%%%%%%%%%%%%%%%%%%%%%%%%%%%%%%%
\section{Ergodicity on Riesz Spaces} \label{ergodicity section}
 
Consider the Riesz space $E=L^1(\Omega,{\cal A},\mu)$, where $\mu$ is a probability measure and $T$ is the expectation operator $T=\E[\cdot]{\bf 1}$, where ${\bf 1}$ is the equivalence class, in $L^1(\Omega,{\cal A},\mu)$, of function with value $1$ almost everywhere. For $\tau$ a measure preserving transformation on $\Omega$ setting $Sf=f\circ \tau$ we have that $(E,T,S,{\bf 1})$ is a conditional expectation preserving system. The definition of ergodicity as $f\circ \tau = f$ for measurable $f$ only if $f$ is almost everywhere constant, can now be written as $Sf=f$ only if $f\in R(T)$. This leads naturally to the following definition in Riesz spaces, for conditional expectation preserving systems.

\begin{definition}[Ergodicity] \label{ergodicity def}
The conditional expectation preserving system $\parent{E,T,S,e}$ is said to be ergodic if $L_S f \in \func{R}{T}$ for all $f \in \3I_S$.
\end{definition}

As will be seen later, this definition preserves the original philosophy of an ergodic process as one in which the time mean (the Ces\`{a}ro mean) is equal to the spatial mean (conditional expectation) in the limit.

We note that $E_e$ is a $R(T)$-module, see \cite{krw}, an hence, when dealing with conditioned systems, the role of the scalars ($\R$ or $\{ke\,|\,k\in\R\}$) is replaced by the ring $R(T)$.

\begin{theorem} \label{tsm thm-1}
The conditional expectation preserving system $\parent{E,T,S,e}$ is ergodic if and only $L_S f=Tf$ for all $f \in \3I_S$.
\end{theorem}

\proof{
Let $\parent{E,T,S,e}$ be ergodic and $f\in \3I_S$, then, by Theorem \ref{bounded riesz birkhoff theorem}, $L_Sf$ exists and $TL_Sf=Tf$. As $\parent{E,T,S,e}$ is ergodic $L_Sf\in R(T)$, so  
there exists $g\in E$ such that $L_Sf=Tg$. As $T$ is a projection $TTg=Tg$, so $Tf=TL_Sf=TTg=Tg=L_Sf.$

Conversely, if $f\in \3I_S$ then, as $\3I_S\subset \3E_S$, $L_Sf$ exists and if we assume $L_Sf=Tf$ then $L_Sf\in R(T)$. Thus $(E,T,S,e)$ is ergodic.
}

For $\parent{E,T,S,e}$ a conditional expectation preserving system we denote by ${\cal B}$ the band projections on $E$, by ${\cal P}:= \setbuilder{P \in \3{B}}{TPe=PTe=Pe}$ the band projections on $E$ which commute with $T$ and by $\3{A}:= \setbuilder{P \in \3{B}}{SPe=Pe}$ the set of $S$ invariant band projections on $E$.
Since $S$ is a positive operator with $Se=e$ and $P$ is positive and below the identity, we have that $SPe, Pe\in E_e$,  for all $P \in \3B$. 
Hence, $L_S$ can always be applied to $SPe$ and $Pe$. We now show that $\3{A}$ characterises $L_S$. 

\begin{theorem}\label{tsm thm}
For $\parent{E,T,S,e}$ a conditional expectation preserving system we have that  ${\cal P}\subset \3{A}=\3{F}$ where $\3{F}:=\setbuilder{P \in \3{B}}{Pe=L_SPe}$ is the set of band projections on $E$ which commute with $L_S$.
\end{theorem}

\proof{
For the purpose of this theorem we can assume that $E=E_e$ as the sets ${\cal A}, {\cal B}, {\cal P}$ do not change when $E$ is replaced by $E_e$. Further,
by Theorem \ref{bounded riesz birkhoff theorem}, $L_S$ restricted to $E_e$ is everywhere defined and a conditional expectation on $E_e$ with $SL_S=L_S$.
 
If $P\in {\cal A}$ then $Pe \in \3I_S$, so, by Lemma \ref{independent exist lemma}, $L_SPe=Pe$, hence $P \in \3F$. Conversely, if $P \in \3F$ then $L_SPe=PL_Se=Pe$, so $SPe=SL_SPe=L_SPe=Pe$, giving $P \in \3A$. Hence, $\3A = \3F$.

From  Theorem \ref{bounded riesz birkhoff theorem},
$L_ST=T=TL_S$ giving that if $P\in{\cal P}$ then $TPe=Pe$ so $L_SPe=L_STPe=TPe=Pe$ and $P\in{\cal F}$.
 }

The following corollary to Theorem \ref{tsm thm} is critical when applying the  concepts of ergodicity in Riesz spaces, further to this, it shows that ergodicity is equivalent to the time and spatial means coinciding in the limit.

\begin{corollary} \label{tsm corollary}
The conditional expectation preserving system $(E,T,S,e)$ is ergodic if and only if $T=L_S$, where $E=\3E_S$.
\end{corollary}

\proof{Theorem \ref{bounded riesz birkhoff theorem} gives that $SL_S=L_S$, hence $R(L_S)\subset \3I_S$.
However, by Lemma \ref{independent exist lemma} $\3I_S\subset R(L_S)$. Hence $\3I_S=R(L_S)$.

If $(E,T,S,e)$ is ergodic, then $Tf=L_Sf=f$, for all $f \in \3I_S$, so $f \in \func{R}{T}$, giving $\func{R}{L_S}=\3I_S \subset \func{R}{T}$. Hence, $\func{R}{L_S}=\func{R}{T}$, but $L_ST=T=TL_S$ so, from \cite{radon paper}, $L_S=T$.

Conversely, suppose that $T=L_S$, then by Theorem \ref{tsm thm-1}, $\parent{E,T,S,e}$ is ergodic.
}

An operator, say $A$, on a Riesz space is said to be strictly positive if $A$ is a positive operator ($Af \geq 0$ for each $f \geq 0$) and $Af=0$ if and only if $f=0$. 
We recall that 
$E_e= \setbuilder{f \in E}{\abs{f} \leq ke \textrm{ for some $k \in {\2R}_+$}},$ 
 the subspace of $E$ of $e$ bounded elements of $E$, 
is an $f$-algebra, see \cite{azouzi, venter, zaanen2}. 
The $f$-algebra structure on $E_e$ is generated by setting $Pe \cdot Qe = PQe$ for all band projections $P$ and $Q$ on $E$ (here, $\cdot$ represents the $f$-algebra multiplication on $E_e$). The linear extension of this multiplication and use of order limits extends this multiplication to $E_e$.
We refer the reader to Azouzi and Trabelsi \cite{azouzi}, Venter and van Eldik \cite{venter} and Zaanen \cite[Chapter 20]{zaanen2} for further background on $f$-algebras.

If $T$ is a conditional expectation operator on $E$ with $Te=e$, then $T$ is also a conditional expectation operator on $E_e$, since, if $f \in E_e$, then $\abs{f} \leq ke$, giving $\abs{Tf} \leq T \abs{f} \leq Tke=ke$. Further $T$ is an averaging operator on $E_e$, that is
$$T(f\cdot g) =f\cdot Tg$$
for $f,g\in E_e$ with $f\in R(T)$, see \cite{expectation paper}.

The following theorem will be seen to be fundamental in linking the concepts of ergodicity and conditional weak mixing in Riesz spaces.

\begin{theorem} \label{thm2}
A conditional expectation preserving system $(E,T,S,e)$, where $T$ being strictly positive, is ergodic if and only if
\begin{equation} \label{thm2 eqn}
 \frac{1}{n} \sum_{k=0}^{n-1} \func{T}{S^kf\cdot g} \to Tf \cdot Tg,
\end{equation}
in order, as $n\to\infty$, for $f,g \in E_e$.
\end{theorem}

\proof{For  $f,g \in E_e$, we note that
	\begin{equation*}
	\frac{1}{n} \sum_{k=0}^{n-1}
	\func{T}{\parent{{S^k}{f}} \cdot g}
	=
	\func{T}{\parent{
	\frac{1}{n} \sum_{k=0}^{n-1} {S^k}{f}
	} \cdot g}
	=
	\func{T}{{S_n}{f} \cdot g},
	\end{equation*}
	and in $E_e$, $S_nf\to L_Sf$, in order, as $n\to\infty$.
	Thus (\ref{thm2 eqn}) is equivalent to 
	\begin{equation}\label{15-sept-1}
	T(L_Sf\cdot g)=Tf\cdot Tg.
	\end{equation}
	
Suppose that $\parent{E,T,S,e}$ is ergodic. As $L_Sf \in \func{R}{T}$ and $T$ is an averaging operator $\func{T}{L_Sf \cdot g}= L_Sf \cdot Tg$. Again from the ergodicity of $\parent{E,T,S,e}$ we have $L_Sf=Tf$.
Thus (\ref{15-sept-1}) holds.

Conversely, suppose that (\ref{thm2 eqn}) or equivalently (\ref{15-sept-1}) holds.
As $T$ is an averaging operator, $Tf\cdot Tg=T(Tf\cdot g)$, which combined with
  (\ref{15-sept-1}) gives 
\begin{equation}\label{15-sept-2}
T((L_Sf-Tf)\cdot g)=0.
\end{equation}
Taking $g=P_\pm e$ where $P_\pm$ are the band projections onto the bands generated by $(L_Sf-Tf)^\pm$ in (\ref{15-sept-2}) gives
\begin{equation}\label{15-sept-2}
0=T((L_Sf-Tf)\cdot P_\pm e)=T(P_\pm (L_Sf-Tf))=T((L_Sf-Tf)^\pm).
\end{equation}
The strict positivity of $T$ now gives that ${\parent{L_Sf-Tf}}^{\pm}=0$, hence $L_Sf=Tf$.
}

%%%%%%%%%%%%%%%%%%%%%%%%%%%%%%%%%%%%%%%%%%%%%%%%%%%%%%%%%%%

\section{Application to conditional weak mixing}

We recall from \cite{khw-1-A} the definition of conditional weak mixing in Riesz spaces. Various other types of mixing in Riesz spaces were studied in \cite{krw}.

\begin{definition}[Weak Mixing] \label{weak mixing definition}
The conditional expectation preserving system  $\parent{E,T,S,e}$ is said to be weakly mixing if,  for all band projections $P$ and $Q$ on $E$,
\begin{equation} \label{weak mixing equation}
 \frac{1}{n} \sum_{k=0}^{n-1} \abs{\func{T}{\parent{S^kPe}\cdot Qe}-TPe \cdot TQe}\to 0,
\end{equation}
in order, as $n\to\infty$.
\end{definition}

From \cite{khw-1-A} we have (\ref{weak mixing equation}) is equivalent to 
\begin{equation}
\label{weak mixing}
\frac{1}{n} \sum_{k=0}^{n-1} \abs{\func{T}{\parent{S^kf}\cdot g}-Tf \cdot Tg} \rightarrow 0
\end{equation}
for all $f,g\in E_e$, in order as $n\to\infty$.

Combining Theorem \ref{thm2} with Lemma \ref{absolute-cesaro} and \eqref{weak mixing}, we have the following.
 
\begin{corollary}
If the conditional expectation preserving system $(E,T,S,e)$, with $T$ strictly positive, is conditionally weak mixing then it is ergodic. 
\end{corollary}

\proof{
If $(E,T,S,e)$ is conditionally weak mixing, then by \eqref{weak mixing},
\[
\frac{1}{n} \sum_{k=0}^{n-1} \abs{\func{T}{\parent{S^kf}\cdot g}-Tf \cdot Tg} \rightarrow 0,
\]
in order, as $n \rightarrow \infty$, for all $f,g \in E_e$, which, by Lemma \ref{absolute-cesaro} gives
\[
\frac{1}{n} \sum_{k=0}^{n-1} \func{T}{\parent{S^kf}\cdot g}-Tf \cdot Tg=\frac{1}{n} \sum_{k=0}^{n-1} \left(\func{T}{\parent{S^kf}\cdot g}-Tf \cdot Tg\right) \rightarrow 0,
\]
in order, as $n \rightarrow \infty$, for all $f,g \in E_e$. That is 
\[
\frac{1}{n} \sum_{k=0}^{n-1} \func{T}{\parent{S^kf}\cdot g}\to Tf \cdot Tg,
\]
in order, as $n \rightarrow \infty$, for all $f,g \in E$, giving that $(E,T,S,e)$ is ergodic, by Theorem \ref{thm2}.
}
%%%%%%%%%%%%%%%%%%%%%%%%%%%%%%%%%%%%%%%

\section{Application to processes with conditionally independent shifts}

We refer the reader to Vardy and Watson, \cite{markov paper}, Kuo, Vardy and Watson, \cite{bernoulli paper}, and Kuo, Labuschagne and Watson, \cite{ergodic paper}, for background on $T$-conditional independence. A summary of $T$-conditional independence, based on \cite{bernoulli paper}, follows.

Let $E$ be a Dedekind complete Riesz space with weak order unit, say, $e$, and $T$ be a conditional expectation operator on $E$ with $Te=e$. Band projections $P$ and $Q$ on $E$ are said to be $T$-conditionally independent (independent with respect to $T$), if
	\begin{equation}
	\label{independence 1}
	TPTQe=TPQe=TQTPe.	
	\end{equation}
We note, from \cite[Lemma 4.2]{grobler-KC}, that band projections $P$ and $Q$ on $E$ are $T$-conditionally independent if and only if 
$TPTQe=TPQe$ (or equivalently $TPQe=TQTPe$), i.e. assuming both inequalities in (\ref{independence 1}) is unnecessary.	
We say that Riesz subspaces $E_1$ and $E_2$ of $E$ are $T$-conditionally independent with respect to $T$ if all band projections $P$ and $Q$ are $T$-conditionally independent where $Pe \in E_1$ and $Qe \in E_2$.
For $\parent{E_{\lambda}}_{\lambda \in \Lambda}$ a family of closed, Dedekind complete Riesz subspaces of $E$ with $\func{R}{T} \subset E_{\lambda}$, for each $\lambda \in \Lambda$, we say that the family is $T$-conditionally independent if, for each pair of disjoint non-empty index sets ${\Lambda}_1, {\Lambda}_2 \subset \Lambda$, we have that $\left\langle \bigcup_{\lambda \in {\Lambda}_1} E_{\lambda} \right\rangle$ and $\left\langle \bigcup_{\lambda \in {\Lambda}_2} E_{\lambda} \right\rangle$ are $T$-conditionally independent. A sequence $\parent{f_n}_{n \in \2N} \subset E$ is said to be $T$-conditionally independent if the family of closed, Dedekind complete Riesz subspaces $\left\langle \braces{f_n} \cup \func{R}{T} \right\rangle$, for each $n \in \2N$, is $T$-conditionally independent.
Here, for $F$ an non-empty subset of $E$, $\left\langle F\right\rangle$ denotes the closed (under order limits in $E$) Riesz subspace of $E$ generated by $F$.
 
We recall the Strong Law of Large Numbers from \cite[Theorem 4.8]{ergodic paper}.
\begin{theorem}[The Strong Law of Large Numbers] \label{strong law of large numbers}
If $\parent{E,T,S,e}$ is a conditional expectation preserving system with $\3E_S=E$ and the sequence $\parent{S^jf}_{j \in \mathbb{N}}$ is independent with respect to $T$, for each $f \in E$, then $T=L_S$.
\end{theorem}

Restricting our attention to $E_e$, with $\displaystyle {\parent{S_n f}}_{n \in \2N}$  independent with respect to $T$, for each $f \in E_e$, we have, from Theorem \ref{strong law of large numbers}, that $T|_{E_e}=L_S|_{E_e}$.
Hence, from Theorem \ref{tsm thm-1} and Theorem \ref{strong law of large numbers}, we have the following corollary.

\begin{corollary}
Let $\parent{E,T,S,e}$ be a conditional expectation preserving system with $\3E_S=E$ and the sequence $\parent{S^jf}_{j \in \mathbb{N}}$  $T$-conditionally independent, for each $f \in E$, then $\parent{E,T,S,e}$ is ergodic.
\end{corollary}

%%%%%%%%%%%%%%%%%%%%%%% Hopf %%%%%%%%%%%%%%%

 \end{document}